\documentclass[reqno]{amsart}

\usepackage{amsmath, amsfonts, amsthm, amssymb, color,  graphicx, mathrsfs, cite}
\usepackage{comment,stmaryrd}

\theoremstyle{plain}
\newtheorem{theorem}{Theorem}[section]
\newtheorem{lemma}{Lemma}[section]
\newtheorem{proposition}{Proposition}[section]
\newtheorem{corollary}{Corollary}[section]

\theoremstyle{definition}
\newtheorem{definition}{Definition}[section]

\theoremstyle{remark}

\newtheorem{example}{Example}[section]

\numberwithin{equation}{section}
\allowdisplaybreaks

\usepackage{ifpdf}
\ifpdf \usepackage[colorlinks=true, citecolor=blue, linkcolor=blue, urlcolor=blue]{hyperref} \fi

\def\thrm{\begin{theorem}}
\def\thrml#1{\begin{theorem}\label{#1}}
\def\ethrm{\end{theorem}}
\def\lmm{\begin{lemma}}
\def\lmml#1{\begin{lemma}\label{#1}}
\def\elmm{\end{lemma}}
\def\dfntn{\begin{definition}}
\def\dfntnl#1{\begin{definition}\label{#1}}
\def\edfntn{\end{definition}}
\def\crllr{\begin{corollary}}
\def\crllrl#1{\begin{corollary}\label{#1}}
\def\ecrllr{\end{corollary}}
\def\xmpl{\begin{example}}
\def\xmpll#1{\begin{example}\label{#1}}
\def\exmpl{\end{example}}
\def\nmrt{\begin{enumerate}}
\def\enmrt{\end{enumerate}}
\def\qtn{\begin{equation}}
\def\qtnl#1{\begin{equation}\label{#1}}
\def\eqtn{\end{equation}}
\def\prpstn{\begin{proposition}}
\def\prpstnl#1{\begin{proposition}\label{#1}}
\def\eprpstn{\end{proposition}}

\def\tm#1{\item[{\rm (#1)}]}
\def\proof{{\bf Proof}.\ }
\def\eprf{\hfill$\square$}

\DeclareMathOperator{\aut}{Aut}

\DeclareMathOperator{\inv}{Inv}

\def\qaq{\quad\text{and}\quad}

\def\lg{\langle}
\def\rg{\rangle}

\def\qaq{\quad\text{and}\quad}

\begin{document}

\title{Notes on finite totally $2$-closed permutation groups}

\author{Gang Chen}
\address{School of Mathematics and Statistics, and Hubei Key Laboratory of Mathematical Sciences, Central China Normal University, P.O. Box 71010, Wuhan 430079, P. R. China.
	}
\email{chengangmath@ccnu.edu.cn}

\author{Qing Ren}
\address{School of Mathematics and Statistics, and Hubei Key Laboratory of Mathematical Sciences, Central China Normal University, P.O. Box 71010, Wuhan 430079, P. R. China.
}
\email{renqing@mails.ccnu.edu.cn}

\thanks{The authors are supported by the NSFC grant No. 11971189.}

\begin{abstract}
Let $N$ be a normal subgroup of a finite group $G$. For a faithful $N$-set $\Delta$, applying the university embedding theorem one can construct a faithful~$G$-set $\Omega$. In this short note, it is proved that if the $2$-closure of $N$ in $\Omega$ is equal to $N$, then the $2$-closure of $N$ in $\Delta$ is also equal to $N$;  in addition, it is proved that any abelian normal subgroup of a finite totally $2$-closed group is cyclic;  finally, it is proved that if a finite nilpotent group is a direct of two nilpotent subgroups where the two factors have coprime orders and both of them are totally $2$-closed then $G$ is totally $2$-closed.  As corollaries, several well-known results on finite totally $2$-closed groups are reproved in more simple ways.  
%\bigbreak
%\noindent {\bf Keywords}: coherent configuration, pseudocyclic association scheme, the Hollmann and Passmann schemes.
\medskip


\medskip

\end{abstract}

%\subjclass[2000]{35B35, 35M20, 35L67, 35Q35.}

\date{}%\today}

\maketitle

%\thispagestyle{empty}
%\vspace{1mm}
%\tableofcontents
%\newpage

\section{Introduction}\label{in}

Let $G$ be a finite group acting faithfully on a finite set $\Omega$. For a positive integer ~$k$, ~$G$ acts naturally on the Catesian product $\Omega^k:=\Omega\times \ldots \times \Omega$. In the year ~1969, Wielandt \cite[Definition 5.3]{W94} introduced the definition of  {\it $k$-closure}~$G^{(k),\Omega}$ of $G\le {\rm Sym}(\Omega)$, which
 is defined as the largest subgroup of ${\rm Sym}(\Omega)$ leaving each orbit of $G$ on $\Omega^k$ invariant. A finite group $G$ is said to be a {\it totally $k$-closed group}, if~$G=G^{(k),\Omega}$ for any faithful $G$-set $\Omega$.

\medskip

In recent years, finite totally $k$-closed groups have been studied in several papers. In \cite{AA16}, it is proved that the center of every finite totally $2$-closed group is cyclic; and finite nilpotent totally $2$-closed groups are characterized. In \cite{CP21}, for a finite abelian group $G$, the minimal positive integer $k$ for which $G$ is totally $k$-closed is given. Actually, the minimal positive integer $k$ is proved to be $1$ plus the number of  invariant factors of $G$. In \cite{A21-1}, it is proved that finite soluble totally ~$2$-closed groups must be nilpotent; the Fitting subgroup of a totally $2$-closed group is shown to be totally $2$-closed. In  \cite{AAPT}, the nontrivial finite totally $2$-closed groups with trivial Fitting subgroups are classified. 

\medskip

In this short note,  we will focus on finite totally $2$-closed groups. Our first result, Theorem A, is a property  on $2$-closures, the second result, Theorem B, is a property on finite totallly $2$-closed goups, and the third result, Theorem C, is on construction of finite nilpotent totally $2$-closed groups. Applying our main results, several results on finite totally $2$-closed groups will be reproved in  self-contained and more simple ways.  

\medskip

Let $G$ be a finite group and $N$ a normal subgroup of $G$. Assume $\Delta$ is a faithful ~$N$-set. By university embedding theorem (ref. Theorem \ref{1129}), we may construct a  faithful~$G$-set $\Omega$. Our first result is the following:

\medskip 

{\bf Theorem A.}\,  {\it Keep the above notation. If $N^{(2), \Omega}=N$, then $N^{(2), \Delta}=N$.}

\medskip

Theorem A captures the common features of the Fitting subgroup, the centralizer of any normal subgroup in a finite totally $2$-closed group. And these subgroups are proved to be totally $2$-closed in a totally $2$-closed group, see Corollary \ref{2007a}. 

\medskip

Base on Theorem A and the known result that every finite totally $2$-closed abelian group is cyclic, we have the second main result of this paper.

\medskip 

{\bf Theorem B.}\, {\it Every normal abelian subgroup of a finite totally  $2$-closed group is cyclic.}

\medskip

Since the classification of finite $p$-groups in which every normal abelian subgroup is cyclic is well known, applying this classification it is proved that a finite nilpotent totally $2$-closed groups is either cyclic or a direct product of a generalized quarternion goup with a cyclic group of odd order. 

\medskip

Additionally, based on the known result that the $2$-closure of any finite nilpotent permutation group is still nilpotent we have the following result, which is the third main result of the paper. 

\medskip 

{\bf Theorem C.}\, {\it Let $G=H\times K$ be a finite nilpotent group with $(|H|, |K|)=1$ and both $H$ and $K$ totally $2$-closed. Then $G$ is totally $2$-closed.} 

\medskip

By Theorem B and Theorem C, the result on the characterization of finite nilpotent totally $2$-closed groups in \cite{AA16} is reproved in a more simple way. 
\medskip

In the final section of this paper, we discuss some known results on $2$-closures of permutation groups in contrast with some results on coherent configurations  . 
 
\medskip

{\bf Notation.} Throughout the paper, $\Omega$ denotes a finite set and the symmetric group of $\Omega$ is denoted by ${\rm Sym}(\Omega)$. The Cartesian product $\Omega\times \Omega$ is written as $\Omega^2$.

 All groups are assumed to be finite nontrivial groups; if a finite group $G$ acts on ~$\Omega$ faithfully, $\Omega$ is also called a faithful $G$-set. 

If a group $G$ acts on $\Omega$, the image of $\alpha\in \Omega$ under the action of $g\in G$ is denoted by $\alpha^g$.

The center of a group $G$ is denoted by $Z(G)$, the centralizer of a subgroup $N$ is denoted by $C_G(N)$, and the Fitting subgroup of $G$ is denoted by $F(G)$.

If $H$ is a subgroup of a group $G$, the core of $H$ in $G$ is dentoed by $H_G$,  i.e., $H_G$ is the intersection of all conjugates of $H$ in $G$.

For a prime $p$, the finite $p$-group $M_{p^{n+1}}$ is defined as 
$$
M_{p^{n+1}}=\lg a,b| a^{p^n}=b^p=1, b^{-1}ab=a^{1+p^{n-1}} \rg.
$$

The cyclic group of order $n$ is denoted by $C_n$.

\medskip

\section{Preliminaries}

Let $G$ be a finite group and $\Omega$ a faithful $G$-set. The {\it $2$-closure} of $G$ is defined as the largest subgroup of ${\rm Sym}(\Omega)$ that leaves every orbit of $G$ on $\Omega^2$ invariant.
\smallskip 
\lmml{1633a}{\rm(\cite[Theorem 5.6]{W94})} Let $\Omega$ be a faithful $G$-set. Then $x\in G^{(2), \Omega}$ if and only if for any $\alpha, \beta\in \Omega$ there exists $g\in G$ such that $(\alpha, \beta)^x=(\alpha, \beta)^g$. 
\elmm

Note that in the above lemma for a fixed element $x\in G^{(2), \Omega}$ the element $g$ depends on the choice of the pair $(\alpha, \beta)\in \Omega^2$.  Also, it is easily seen that $G^{(2), \Omega}$ contains $G$ as a subgroup.

\lmml{1919c}{\rm(\cite[Lemma 2.2]{AA16})}Let $G\le {\rm Sym}(\Omega)$ and $H\le {\rm Sym}(\Gamma)$, and let $(\theta, \lambda)$ be a permutation isomorphism from  $(G, \Omega)$ to $(H, \Gamma)$. Then $(\theta, \lambda)$ can be extended a permutation isomorphism $(\tilde{\theta}, \tilde{\lambda})$ from $(G^{(2), \Omega}, \Omega)$ to $(H^{(2), \Gamma}, \Gamma)$, where $\tilde{\lambda}=\lambda$ and the restriction of $\tilde{\theta}$ to $G$ is equal to $\theta$.  
\elmm

\lmml{1946c}{\rm(\cite[Lemma 2.1]{AA16}, \cite[Lemma 1.1]{A21-1})} Let $G$ be a finite totally $2$-closed group and  $N$  a subgroup of $G$. Then $C_G(N)^{(2), \Omega}=C_G(N),$ for any finite faithful ~$G$-set $\Omega$.
\elmm

The following lemma is a special case of Lemma 3.2 of \cite{CP21}. For completeness, we give a proof. 

\lmml{908c}{\rm (\cite[Lemma 3.2]{CP21})}Let $m,n$ be positive integers with  $m$ divisible by $n$. Then $C_n\times C_m$ is not totally $2$-closed. 
	\elmm
	\proof Let $\Omega=\{1,2, \ldots, 2n+m\}$. Choose two permutations in ${\rm Sym}(\Omega)$ in the following way:
	$$
	x_1=(1\ldots n)(n+1,\ldots,2n), x_2=(n,n-1,\ldots, 1)(2n+1,\ldots, 2n+m).
	$$
	One can easily see that the subgroup $H=\lg x_1,x_2\rg$ of ${\rm Sym}(\Omega)$ is isomorphic to the direct product $C_n\times C_m$. However, the permutation $(1\ldots n)$ belongs to $H^{(2), \Omega}$ but not to $H$. The lemma then follows. 
	
	\eprf
	
	\lmml{1107a}{\rm(\cite[Lemma 2.1]{A21-1})}
	Let $G=HK$ be a totally 2-closed group, where ~$H$ and $K$ are proper subgroups of $G$. If $H_G\cap K_G=1$, then $G=H_G\times K_G$. Furthermore,  $H_G$ and $K_G$ are totally 2-closed. In particular, if $H\cap K=1$, then  $G=H\times K$ and both $H$ and $K$ are totally 2-closed.  
	\elmm 
	
\smallskip	
	
\lmml{952b}{\rm (\cite[Theorem 3]{AA16})}If $G$ is a finite totally $2$-closed abelian group, then ~$G$ is cyclic. 
\elmm
\proof It suffices to show that for every prime divisor $p$ of $|G|$, the Sylow $p$-subgroup $S_p$ of $G$ is cyclic. 

Suppose on the contrary that the Sylow $p$-subgroup~ $S_p$ of $G$ is not cyclic for some prime divisor $p$ of $|G|$. Then there exist positive integers ~$e_1\ge e_2$ and a subgroup $K$ of $G$ such that 
$$
G=(C_{p^{e_1}}\times C_{p^{e_2}})\times K.
$$
By Lemma \ref{1107a}, $C_{p^{e_1}}\times C_{p^{e_2}}$ must be totally $2$-closed. This is a contradiction to Lemma \ref{908c}. 
\eprf

\smallskip

\lmml{1110b}{\rm (\cite[Chap. 3, Theorem 7.6]{Hu67})}Let $p$ be a prime and $G$ a finite nontrivial~ $p$-group such that every normal abelian subgroup of $G$ is cyclic. Then 
\nmrt
\tm{i} If $p>2$, then $G$ is cyclic. 
\tm{ii} If $p=2$, then $G$ has a cyclic subgroup of index $2$.
\enmrt
\elmm

\smallskip

The following lemma is based on the fact that as permutation groups cyclic groups of prime power order and generalized quaternion groups have trivial one point stabilizers.

\lmml{1018a}{\rm(\cite[Lemma 2.5]{AA16})} Finite cyclic groups of prime power orders and generalized quaternion groups are totally $2$-closed. 
\elmm

\lmml{1511c}{\rm(\cite[Corollary 3.4]{A21-1})} Let $G\le {\rm Sym}(\Omega)$. Then $G$ is nilpotent if and only if $G^{(2),\Omega}$ is nilpotent. 
\elmm

\section{Main Results}

The following result plays a key role in the whole theory. 

\smallskip

\thrml{1129}{\rm (University embedding theorem \cite[Theorem 2.6 A]{DM96})} Let $G$ be an arbitrary group with a normal subgroup $N$ and put $K:=G/N$. Let $\psi: G\rightarrow K$ be a homomorphism of $G$ onto $K$ with kernel $N$. Let $T:=\{t_u|u\in K \}$ be a set of right coset representatives of $N$in $G$ such that $\psi(t_u)=u$ for each $u\in K$. Let $x\in G$ and~ $f_x:\, K\rightarrow N$ be the map with $f_x(u)=t_uxt^{-1}_{u\psi(x)}$ for all $u\in K$. Then~ $\phi(x)=: (f_x, \psi(x))$ defines an embedding $\phi$ of $G$ into $N\wr K$. Furthermore, if~ $N$ acts faithfully on a set ~$\Delta$ then $G$ acts faithfully on $\Omega:=\Delta\times K$ by the following rule~ $(\delta, k)^x=(\delta^{f_x(k)}, k\psi(x))$. 
\ethrm

In the first paragraph on page 9 of \cite{A21-1}, it was shown that the action of $N$ on $\Omega$ is permutation isomorphic to the natural action of $N$on a disjoint union of $|G:N|$ copies of $\Delta$.  Here, for each $g\in K$,  the corresponding copy may be denoted as 
$$
\Delta_g:=\{(\delta, g):\delta\in \Delta\}, 
$$
and the permutation isomorphism is given by  
$$
(\theta_g, \lambda_g): \, (N, \Delta_g)\, \rightarrow \, (N, \Delta), \quad n\mapsto t_gnt_g^{-1}, \, (\delta, g)\mapsto \delta.
$$
Then $N$ acts on $\Omega$ in the following way: for any $\delta\in \Delta$, $g\in K$, and $x\in N$, 
$$
(\delta, g)^x=(\delta^{\theta_g(x)}, g).   
$$

One can obtain a permutation isomorphism  $(\mu_g, \lambda_g): (N^{(2), \Delta_g}, \Delta_g) \rightarrow (N^{(2), \Delta}, \Delta)$, by Lemma \ref{1919c},  for each $g\in K$. Also, the restriction of $\mu_g$ to $N$ is equal to $\theta_g$. Furthermore, ~$N^{(2), \Delta_g}$ acts faithfully on $\Omega$ in the following way: 
$$
(\delta, g)^x=(\delta^{\mu_g(x)}, g),\quad \forall \delta\in \Delta, g\in K, \qaq x\in N^{(2), \Delta_g} .
$$

For any $x\in N^{(2), \Delta}$, denote the images of $x$ in ${\rm Sym}(\Omega)$ by $\tilde{x}$. For any $(\delta, g)\in \Omega$, 
\qtnl{947b}
(\delta, g)^{\tilde{x}}=(\delta, g)^{\mu_g^{-1}(x)}=(\delta^x,g).
\eqtn

\smallskip

The following theorem, which is Theorem A in this paper,  is a generalization of \cite[Proposition 3.6]{A21-1}. 

\thrml{1000a} In the context of Theorem \ref{1129} , if we keep the above notation and assume that  $N^{(2), \Omega}=N$, then $N^{(2), \Delta}=N$. 
\ethrm 

\proof It suffices to show that, as permutation groups on $\Omega$,    the cardinality of the image of $N^{(2), \Delta}$ in ${\rm Sym}(\Omega)$ is not greater than  that of $N$. 
 
\smallskip

For any $x\in N^{(2), \Delta}$, and any $(\delta_1, g_1), (\delta_2,g_2)\in \Omega$, there exists $n\in N$ such that $(\delta_1, \delta_2)^x=(\delta_1, \delta_2)^n$. Denote the images of $x$ and $n$ in ${\rm Sym}(\Omega)$ by $\tilde{x}$ and $\tilde{n}$, respectively. 
By formula \eqref{947b}, one can see that  
\begin{equation*}
\begin{split}
((\delta_1, g_1), (\delta_2, g_2))^{\tilde x} &=((\delta_1^{x}, g_1), (\delta_2^{x}, g_2))\\
&=((\delta_1^{n}, g_1), (\delta_2^{n}, g_2))\\
&=((\delta_1, g_1), (\delta_2, g_2))^{\tilde n}.
\end{split}
\end{equation*}
We conclude that ${\tilde x}\in N^{(2), \Omega}=N$. Thus, the cardinality of the image of $N^{(2), \Omega}$ is not greater than that of $N$, as required. 

\eprf

\crllrl{2007a}{\rm (\cite[Theorem B, Proposition 3.7]{A21-1})} Let $G$ be a finite totally $2$-closed group and $N$ a normal subgroup of $G$. Then $C_G(N)$, $Z(G)$, and $F(G)$ are totally $2$-closed group. 
\ecrllr
\proof The corollary is a consequence of Lemma \ref{1946c}, Theorem \ref{1000a}, and \cite[Lemma 3.5]{A21-1}.

\eprf

\lmml{1014a}{\rm (\cite[Theorem 1]{AA16})}The center of every totally $2$-closed group is cyclic. 
\elmm
\proof By Corollary \ref{2007a}, the center is totally $2$-closed and thus is cyclic by Lemma \ref{952b}. 
\eprf

\smallskip

The following result is the second main result, Theorem B,  in this paper. 

\smallskip

\thrml{1102b}Let $G$ be a finite totally $2$-closed group and $N$ a normal abelian subgroup of $G$. Then $N$ is cyclic. 
\ethrm
\proof By Corollary \ref{2007a}, $C_G(N)$ is totally $2$-closed. Thus the center of $C_G(N)$ is cyclic by Lemma \ref{1014a}. Note that $N$ is a central subgroup of $C_G(N)$. It then follows that $N$ is cyclic. 
\eprf

\thrml{1135c}Let $p$ be a prime and $G$ is a finite totally $2$-closed nontrivial $p$-group. Then either $G$ is cyclic or $G$ is a generalized quaternion group. 
\ethrm
\proof Assume that $G$ is a finite totally $2$-closed $p$-group and that $G$ is not cyclic.

By Theorem \ref{1102b}, every normal abelian subgroup of $G$ is cyclic.  By Lemma \ref{1110b}, we see that $p=2$ and $G$ has a cyclic subgroup of index $2$. By \cite[Theorem 1.2]{B08}, $G$ is either a dihedral group, a generalized quaternion group, ~$M_{2^{n+1}}$, or a semidihedral group. 

However, except for the generalized quaternion group, the other three groups are semidirect products of two proper subgroups. By Lemma \ref{1107a}, if $G$ is equal to a semidirect product of two proper subgroups, then $G$ must be the direct product of these two proper subgroups, which is a contradiction.  We thus conclude that $G$ is a generalized quaternion group.  

\eprf

As a consequence of  Lemma \ref{1107a} and Theorem \ref{1135c}, we obtain the following.

\crllrl{1041b} If $G$ is a finite nilpotent totally $2$-closed group, then $G$ is cyclic or a direct product of a generalized quaternion group with a cyclic group of odd order. 
\ecrllr

The following theorem is the third main result, Theorem C, in this paper. 

\thrml{1519b}Let $G=H\times K$ be a finite nilpotent group with $(|H|, |K|)=1$ and both $H$ and $K$ totally $2$-closed. Then $G$ is totally $2$-closed. 
\ethrm
\proof Let $\Omega$ be a faithful $G$-set. By Lemma \ref{1511c}, $G^{(2), \Omega}$ is nilpotent. We may write
$$
G^{(2), \Omega}=J\times L, 
$$
where $(|J|, |L|)=1$,  $H\le J$, and $K\le L$.   

Set $j:=|J|$ and $l:=|L|$. For any $x\in G^{(2),\Omega}$ and any $(\alpha, \beta)\in \Omega^2$, there exists ~$g\in G$ such that $(\alpha, \beta)^x=(\alpha, \beta)^g$ by Lemma \ref{1633a}.  By the assumption, one can write $x=yz$ with $y\in J$ and $z\in L$, and $g=hk$ with $h\in H$ and $k\in K$. Thus, 
$$
(\alpha, \beta)^{x^l}=(\alpha, \beta)^{g^l}\, \Rightarrow \, (\alpha, \beta)^{y^l}=(\alpha,\beta)^{h^l}.
$$ 
This implies that $y^l\in H^{(2), \Omega}$ by Lemma \ref{1633a}. As we are assuming that $H$ is totally ~$2$-closed, $y^l\in H$. In addition, the fact that $(l, |H|)=1$ yields that $y\in H$. Similarly, one can show that $z\in K$. 

\smallskip

It follows that $x=yz\in HK=G$ and therefore $G^{(2), \Omega}=G$. Since this is true for any faithful $G$-set, $G$ is totally $2$-closed. 

\eprf

\crllrl{940a} Let $G$ be a finite cyclic group or a direct of a generalized quaternion group with a cyclic group of odd order. Then $G$ is totally $2$-closed. 
\ecrllr
\proof By the assumption we may decompose $G$ into a direct product of its Sylow subgroups, where each Sylow subgroup is totally $2$-closed by Lemma \ref{1110b}. Thus, $G$ is totally $2$-closed by Theorem \ref{1519b}.  
\eprf

\smallskip

As consequences of Corollaries \ref{1041b} and \ref{940a}, we have the following. 

\crllrl{955b}{\rm (\cite[Theorem 2]{AA16})} A finite nilpotent group is totally $2$-closed if and only if it is cyclic or a direct product of a generalized quaternion group with a cyclic group of odd order. 

\ecrllr

\section{Concluding Remarks}

Let $G\le {\rm Sym}(\Omega)$. Then the action of $G$ on $\Omega^2$ generates a combinatorial structure, i.e., the {\it coherent configuration} associated with $G$, denoted by $\inv (G)$. The ~$2$-closure of $(G, \Omega)$ is then the automorphism group of the coherent configuration. In other words, $G^{(2), \Omega}=\aut(\inv(G))$, see \cite[Definition 2.2.14]{CP}. 

\smallskip

Several known results on $2$-closures of finite permutation groups have their correspondents in the theory of coherent configurations. For example,  Theorem 3 in \cite{A21-1} corresponds to Theorem 3.2.21 and Theorem 3.4.6 in \cite{CP}.  And Corollary 3.2 in  \cite{A21-1} corresponds to Exercise 2.7.17 $(4)$ in \cite{CP}. Finally, Lemma 1.4 in \cite{A21-1} or Lemma 2.9 in \cite{AA16} is a consequence of Theorem 3.2.5 in \cite{CP}. More precisely, if $G_i\le {\rm Sym}(\Omega_i)$ and ~$\Omega$ is the disjoint union of $\Omega_i's$, $i=1, \ldots, n$, then $G=G_1\times\cdots \times G_n$ acts  faithfully on $\Omega$ and $G^{(2), \Omega}=G_1^{(2), \Omega_1}\times \cdots \times G_n^{(2), \Omega_n}$ by the cited Theorem.

\medskip

{\bf Acknowledgements.} \, The authors would like to thank the refree for valuable comments.


\begin{thebibliography}{99}\small
	
\bibitem{AA16} A. Abdollahi and M.Arezoomand, \emph{Finite nilpotent groups that coincide with their $2$-closures in all of their faithful permutation representations}. J. Algebra Appl. {\bf 17(4)} (2018), 1850065.

\bibitem{AAPT}
M.Arezoomand, M. A. Iranmanesh,  C. E. Praeger and G.Tracey, \emph{Totally $2$-closed finite groups with trivial Fitting subgroup}, arxiv preprint, arxiv:2111.02253, 2021.

\bibitem{A21-1}
A.~ Abdollahi, M.~Arezoomand, and G.Tracey,  \emph{On finite totally $2$-closed groups}, arxiv preprint, arxiv:2001.09597, 2021.

\bibitem{B08}
Y. Berkovich, \emph{Groups of prime power order,} Vol. 1, Expositions in Mathematics 46, Walter de Gruyter GmbH $\&$ Co.KG, Berlin, 2008.  

\bibitem{CP21}
D. Churikov and C. E. Praeger,  \emph{Finite totally $k$-closed groups}, Tr. Inst. Mat. Mekh. {\bf 27 (1)}(2021), 240-246.


\bibitem{CP} 
G.~Chen and I.~Ponomarenko, \emph{Coherent Configurations}, Central China Normal Universit Press, Wuhan, 2019.

\bibitem{DM96}
J.D. Dixon and B. Mortimer, \emph{Permutation groups}, Springer, New York, 1996.

\bibitem{Hu67} 
B. Huppert, \emph{Endliche Gruppen. I.} Die Grundlehren der Mathematischen Wissenschaften, Band 134. Spring-Verlag, Berlin-New York, 1967.
    
   
\bibitem{W64} H. Wielandt, \emph{Finite permutation groups,} Academic Press, New York, 1964.


\bibitem{W94}
 H.W. Wielandt, \emph{`Permutation groups through invariant relations and invariant functions'}, Lecture Notes, Ohio State University, 1969. Also published in: Wielandt, Helmut, Mathematische Werke (Mathematical works) Vol. 1. Group theory, Walter der Gruyter $\&$  Co., Berlin, 1994, pp. 237-296. 


\end{thebibliography}
\end{document}